\documentclass[11pt]{amsart}
\usepackage{amssymb}
\usepackage{amsmath}
\usepackage{amsfonts,amsthm,latexsym,mathrsfs,bbm}
\usepackage[usenames,dvipsnames,svgnames,table]{xcolor}
\usepackage{hyperref}
\hypersetup{colorlinks=true,linkcolor=red,citecolor=gray}
\usepackage{amsthm}
\theoremstyle{plain}

\newtheorem*{theorem*}{{\bf Theorem}}
\newtheorem*{maintheorem*}{{\bf Main Theorem}}
\newtheorem*{lemma*}{{\bf Lemma}}
\newtheorem{corollary}[subsubsection]{{\bf Corollary}}
\newtheorem{proposition}[subsubsection]{{\bf Proposition}}
\newtheorem{lemma}[subsubsection]{{\bf Lemma}}
\newtheorem*{claim}{{\bf Claim}}
\newtheorem*{remark*}{{\bf Remark}}

\theoremstyle{remark}

\newtheorem{example}[subsubsection]{{\it Example}}
\numberwithin{equation}{subsection}

\DeclareMathOperator{\ad}{ad}
\DeclareMathOperator{\im}{im}

\DeclareMathOperator{\Aut}{Aut}
\DeclareMathOperator{\Fix}{Fix}
\DeclareMathOperator{\End}{End}

\def \FF {\mathbb F}
\def \Prob {\mathbb P}

\DeclareMathOperator{\Alt}{Alt}

\DeclareMathOperator{\PSL}{PSL}
\DeclareMathOperator{\SL}{SL}

\DeclareMathOperator{\Sz}{Sz}

\newcounter{ithmcount}

\makeatletter
\def\@author#1{\g@addto@macro\elsauthors{\normalsize%
    \def\baselinestretch{1}%
    \upshape\authorsep#1\unskip\textsuperscript{%
      \ifx\@fnmark\@empty\else\unskip\sep\@fnmark\let\sep=,\fi
      \ifx\@corref\@empty\else\unskip\sep\@corref\let\sep=,\fi
      }%
    \def\authorsep{\unskip,\space}%
    \global\let\@fnmark\@empty
    \global\let\@corref\@empty
    \global\let\sep\@empty}%
    \@eadauthor={#1}
}
\makeatother

\title[Gaps in word probabilities]
{Gaps in probabilities of satisfying some commutator-like identities}

\author[C. Delizia]{Costantino Delizia}
\address{University of Salerno, Italy}
\email{cdelizia@unisa.it}

\author[U. Jezernik]{Urban Jezernik}
\address{University of the Basque Country, Spain}
\email{urban.jezernik@ehu.eus}

\author[P. Moravec]{Primo\v z Moravec}
\address{University of Ljubljana, Slovenia}
\email{primoz.moravec@fmf.uni-lj.si}

\author[C. Nicotera]{Chiara Nicotera}
\address{University of Salerno, Italy}
\email{cnicoter@unisa.it}

\date{\today}

\thanks{
The first and the fourth author have been supported by the \lq\lq National Group for Algebraic and Geometric Structures, and their Applications\rq\rq\ (GNSAGA -- INdAM).
The second author has received funding from the European Union’s Horizon 2020 research and innovation programme under the Marie Sklodowska-Curie grant agreement No 748129. He has also been supported by the Spanish Government grant MTM2017-86802-P and by the Basque Government grant IT974-16.
The third author has been partially supported by the Slovenian Research Agency (research core funding No. P1-0222, and projects No. J1-8132, J1-7256 and N1-0061).}

\begin{document}
\baselineskip=15pt

\begin{abstract}

We show that there is a positive constant $\delta < 1$ such that the
probability of satisfying either the $2$-Engel identity $[X_1, X_2, X_2] = 1$ or
the metabelian identity $[[X_1, X_2], [X_3, X_4]] = 1$ in a finite group
is either $1$ or at most $\delta$.

\end{abstract}

\maketitle



\section{Introduction} \label{s:intro}


It is an old, elegant, well-known, and at the same time somewhat surprising
result that the probability that two randomly chosen elements commute in a
nonabelian finite group can not be arbitrarily close to $1$. To be more
precise (see \cite{Gus73}), the commuting probability of no finite group can
belong to the interval $(\frac58, 1)$, and so there is a gap in the possible
probability values. Following on this, many other, deeper results on the
structure of the set of all possible values of the probability of satisfying
the commutator identity have since emerged (see \cite{Ebe15} and the references
therein).


Recently, more general word maps on finite groups have been explored from a standpoint of a similar probabilistic flavour. Here, a word map on a group $G$ is a map $w
\colon G^d \to G$ induced by substitution from a word $w \in F_d$ belonging to
a free group of rank $d$.  For a fixed element $g \in G$ of a finite group $G$, set
\[
\Prob_{w=g}(G) =
\frac{|w^{-1}(g)|}{|G|^d}
\]
to be the probability that $w(g_1, g_2, \dots, g_d) = g$ in $G$, where $g_1,
g_2, \dots, g_d$ are chosen independently according to the
uniform probability distribution on $G$.
Following recent breakthroughs on the values of these probabilities for finite
simple groups (see \cite{LarSha12}, \cite{Nik16} and \cite{Bor17}), applications have been
developed also for infinite groups
(see \cite{LarSha17} for an approach via the Hausdorff dimension for
residually finite groups),
indicating how these probabilities of finite quotients of
a given infinite group are tightly related to its algebraic structure
(see \cite{Sha18}, as well as \cite{MTVV18} for a more
geometric approach).
This has been done mostly for the simple and longer commutator words,
the results ultimately resting on the aforementioned probability gap
and its stronger variants (see \cite{Neu89}).


The purpose of this paper is to investigate some words that are natural
generalizations of the simple commutator. 
Our main
result shows the existence of gaps in probabilities of satisfying these words.

\begin{maintheorem*}
Let $w$ be either the $2$-Engel or the metabelian word.
There exists a constant $\delta < 1$ such that whenever $w$
is not an identity in a finite group $G$, we have $\Prob_{w=1}(G) \leq \delta$.
\end{maintheorem*}

The strategy of the proof is quite general and might be applied to some other
words. 
In particular, it applies to the long commutator word and thus provides
an alternative proof to its already known probability gap (see \cite{Erf07}).
The general strategy goes as follows.
Suppose $G$ is a finite group in which $w$ is not
an identity. Consider a chief factor $G_1/G_2 = T^k$ of $G$, where $T$ is a
simple group and $k \geq 1$.

If $T$ can be chosen to be nonabelian, then $T$ does not satisfy the word $w$.
For the purposes of our claim, we can replace $G$ by its quotient
$G/C_G(G_1/G_2)$ and therefore assume that $T^k \leq G
\leq \Aut(T^k)$ (see Section \ref{s:nonsolvable_groups}). 
Word probabilities in such groups have been studied
extensively. As long as $T$ is a large enough simple group, our claim follows
from the works of Larsen and Shalev (see \cite{LarSha17}). For small groups $T$, it
is necessary to bound their multiplicities $k$ in $G$. This is tightly related
to the concept of varied coset identities of $T$ (see \cite{Nik16}).  By
inspecting the required cases following Bors (see \cite{Bor17}),
we are able to achieve the
goal of bounding these multiplicities
for both the $2$-Engel word and the metabelian word (see Subsection
\ref{subs:multiplicity_bounding_words}).

\begin{theorem*}
The $2$-Engel and the metabelian word are
multiplicity bounding.
\end{theorem*}

The proof relies of inspecting fixed points of outer automorphisms
of finite simple groups for the $2$-Engel word, whereas a more direct
algebraic manipulation works for the metabelian word.
Along the way, we also investigate coset probabilities of finite groups
with respect to nonsolvable normal subgroups
(see Subsection \ref{subs:coset_probabilities}).
We show that in infinite
groups that possess infinitely many nonabelian upper composition factors,
the infimum of probabilities of satisfying the $2$-Engel or metabelian word
in its finite quotients is $0$.

On the other hand, if all the chief factors of $G$ are abelian, then $G$ is
solvable. In this case, the verbal subgroup $w(G)$ can be assumed to be the
unique minimal normal subgroup of $G$, and so it is a vector space over a
finite field. We proceed by analysing the linear representation of $G$ on this
space. As long as this representation is nontrivial, we are able to provide a
general procedure on how to establish a word probability gap
(see Subsection
\ref{subs:general_principle_for_bounding_probability_higher_dimensional_verbal}).
This is then
executed for the $2$-Engel word, where the only problematic elements are those
acting quadratically on $w(G)$, and for the metabelian word, where the
situation is simpler due to invariance of variables. In the case when the
representation of $G$ on the verbal subgroup is $1$-dimensional or, in the case of the metabelian word, can be factorized through the abelianization of $G$, 
it is not possible to obtain any
information from the representation alone. Here, we instead consider the
restriction of the word map on the coordinate axes
(see Subsection
\ref{subs:general_principle_for_bound_probability_1_dimension}).
An argument involving the
analysis of whether or not such an induced map is trivial works for both the
$2$-Engel word and the metabelian one. Joint with the above Theorem, we
conclude the validity of the Main Theorem.

An explicit value of $\delta$ in the Main Theorem could, in
principle, be determined by examining the proof. For this, one would need to
compute the probabilities $\Prob _{w=1}(G)$ for all finite groups $G$ with
$T^k \leq G \leq \Aut(T^k)$, where $T$ is a nonabelian finite simple group and
both $|T|$ and $k$ are bounded in terms of $w$. The difficulty lies in finding
good bounds. On the other hand, when restricting only to solvable groups, it
follows from our proofs that, for the $2$-Engel word, one can take $\delta =
\frac34$, equality being attained with the dihedral group $D_{16}$, and for
the metabelian word, one has $\delta \leq \frac{29}{32}$, but this
bound might not be sharp.

\smallskip

A word on the notation. The generators of the free group $F_d$ will be denoted
by $X_1, \dots, X_d$. The multiplicity of $X_i$ in $w \in F_d$ will be denoted by
$\mu_w(X_i)$. The length of $w$ will be denoted by $\ell(w)$.

\section{Nonsolvable groups}
\label{s:nonsolvable_groups}

In this section, we deal with bounding the probability $\Prob_{w=1}(G)$ for
finite nonsolvable groups $G$. We will repeatedly use the following reduction
lemma.

\begin{lemma*}
Let $w$ be a nontrivial word. Let $G$ be a finite group and $N$ a normal
subgroup of $G$. Then
$\Prob_{w=g}(G) \leq \Prob_{w=gN}(G/N)$
for every $g\in G$.
\end{lemma*}

We can therefore replace $G$ by its quotient
$G / C_G(T^k)$ and hence reduce our claim to bounding
the probability
$\Prob_{w=1}(G)$ in the case when $T^k \leq G \leq \Aut(T^k)$ for a nonabelian
finite simple group $T$.

\subsection{Large simple groups}

\begin{lemma}[\cite{LarSha17}, Theorem 1.8]
	\label{l:larsen}
    Let $G$ be a finite group such that $T^k \le G \le  \Aut (T^k )$ for some $k
\ge 1$ and a finite nonabelian simple group $T$. Suppose $w$ is a nontrivial
word. Then there exist constants $C = C (w)$, $\epsilon = \epsilon (w) > 0$
depending only on $w$ such that, if $|T| \geq C$, then for any $g \in G$ we have
$\Prob_{w=g}(G) \leq |T^k|^{-\epsilon}$.
\end{lemma}

As long as $|T| > C$, we therefore have $\Prob_{w=1}(G) < C^{- \epsilon}$.

\subsection{Multiplicity bounding words}
\label{subs:multiplicity_bounding_words}

A reduced word $w$ is called {\em multiplicity bounding} (see \cite{Bor17}) if,
whenever $G$ is a finite group such that $\Prob_{w=g}(G) > \rho$ for some $g \in G$,
the multiplicity of a nonabelian simple group $S$ as a composition factor
of $G$ can be bounded above by a function of only $\rho$ and $S$.

Whenever our word $w$ is multiplicity bounding, we can solve our problem for
groups $G$ with $|T| \leq C$. Namely, for each of these nonabelian groups $T$,
we either have that  $\Prob_{w = 1}(G) \leq \frac12$, in which case we are done, or
we can assume that $\Prob_{w = 1}(G) > 1/2$. In the latter case, the multiplicity
$k$ of $T$ is bounded above by a constant depending only on $w$ and $T$.  This
means that $|G| \leq \lvert \Aut(T^k) \rvert$ is bounded above by a constant, and we
therefore have an upper bound for $\Prob_{w=1}(G)$ as well.
Note that since $w$ is multiplicity bounding, it can not be
an identity in any of these finitely many groups $G$
(since it would otherwise also be an identity in some nonabelian
finite simple group $S$ and thus in all powers $S^n$ as well).

We have thus proved the following.

\begin{proposition} \label{p:gap_nonsolvable}

Let $w \in F_d$ be a multiplicity bounding word.
Then there exists a constant $\delta = \delta(w) < 1$ such
that every nonsolvable finite group $G$ satisfies $\Prob_{w = 1}(G) \leq \delta$.

\end{proposition}

It therefore remains to deal with proving that the words we are interested in
are indeed multiplicity bounding. In fact, we will prove that these words
satisfy a stronger property.

\subsection{Coset word maps and variations}

Let $S$ be a nonabelian finite simple group. The word $w \in F_d$ defines a {\em
word map} $S^d \to S$ by evaluation $(s_1, \dots, s_d) \mapsto w(s_1, \dots,
s_d)$. Consider $S \leq \Aut(S)$ and let $g_1, \dots, g_d \in \Aut(S)$.
Then there is a corresponding {\em coset word map} $S^d \to S$ defined by
$(s_1, \dots, s_d) \mapsto w(s_1 g_1, \dots, s_d g_d)$.

We will also require the notion of a {\em variation} of $w$. This is a word
$\tilde w$ obtained from $w$ by adding, for each $1 \leq i \leq d$, to each
occurrence of $X_i^{\pm 1}$ in $w$ a second index from the range $\{ 1, 2, \dots,
\mu_w(X_i) \}$. To each such variation $\tilde w$, we can associate a {\em
varied coset word map} of $w$, which is just a coset word map of the variation
$\tilde w$.

\subsection{Very strongly multiplicity bounding words}

A word $w$ is called {\em very strongly multiplicity bounding} (VSMB)
if for all nonabelian finite simple groups $S$, none of the varied coset
word maps of $w$ on $S$ is constant. Such words are multiplicity bounding
(see \cite[Proposition 2.9]{Bor17}).

\begin{example}[\cite{Bor17}, Corollary 3.4]

Long commutator words $\gamma_n(X_1, \dots, X_n) = [X_1, \dots, X_n]$ are
all VSMB.

\end{example}

The following criterion for being VSMB will be of use.

\begin{lemma}[\cite{Bor17}, Proposition 3.1, Proposition 6.1]
\label{l:bors_reductions_for_VSMB}

The following words are VSMB.

\begin{enumerate}
\item Words in which some variable occurs with multiplicity $1$.
\item Words in which some variable occurs with multiplicity $2$, provided that
either $w = w_1 X_d w_2 X_d w_3$
or $w = w_1 X_d^{\pm 1} \bar w_2 X_d^{\mp 1} w_3$,
where $w_1, w_2, w_3$ are reduced and $\bar w_2$ is VSMB.
\item Words of length at most $8$ excluding the power word $X_1^8$.
\end{enumerate}

\end{lemma}

In order to verify that a given word is VSMB, it suffices to inspect only a
limited set of simple groups. The following criterion will suffice here.

\begin{lemma}[\cite{Bor17}, Proposition 4.9 (8)]
\label{l:bor_reduction_simple_groups}

Let $w \in F_d$ be a reduced word. Set $m = \max_i \mu_w(X_i)$. Then $w$ is VSMB
as long as the word map of $w$ on $\PSL_2(2)$ and $\Sz(2)$ is not constant and
none of the varied coset words maps on $\PSL_2(p^n)$ for a prime 
$2 < p \leq m$ and $n$ a power of $2$ is constant.



\end{lemma}

\subsection{Automorphisms of the relevant simple groups}

We will be inspecting coset word maps on the simple groups from Lemma
\ref{l:bor_reduction_simple_groups}. For this, we will need to understand
cosets of inner automorphisms of these groups.

The automorphisms of $\PSL_2(2^n)$ and $\Sz(2^n)$ consist of inner automorphisms
and field automorphisms. The field automorphisms are generated by the Frobenius
automorphism $\sigma$ that extends the field automorphism  $\FF_{2^n} \to
\FF_{2^n}$, $x \mapsto x^2$.  The order of $\sigma$ in $\Aut(S)$ is equal to
$n$, which is assumed to be a prime.

As for the groups $\PSL_2(p^n)$ with $p$ odd, there is an additional
outer automorphism $D \in \Aut(S)$ induced by conjugation with the diagonal
matrix
\[
\begin{pmatrix}
\omega & 0 \\
0 & 1
\end{pmatrix},
\]
where $\omega$ is a generator of $\FF_{2^n}^\times$. This automorphism
fixes all the diagonal matrices. Moreover, it
satisfies the relation $[D, \sigma] = D^{p-1}$.  Therefore every element of the
group $\langle D, \sigma \rangle$ can be written uniquely as $\sigma^i D^j$ with
$0 \leq i < n$ and $0 \leq j < p^n - 1$. The square $D^2$ is  an inner
automorphism, given as conjugation with  the matrix
\[
\begin{pmatrix}
\omega^2 & 0 \\
0 & 1
\end{pmatrix} \equiv
\begin{pmatrix}
\omega & 0 \\
0 & \omega^{-1}
\end{pmatrix} \in \PSL_2(p^n).
\]

We will require the following property of these outer automorphisms.

\begin{lemma} \label{l:out_fixed_points}
Let $1 \neq \alpha \in \langle D, \sigma \rangle \leq \Aut(\PSL_2(p^n))$
with $p > 2$ and $n$ a power of $2$. Then
\[
\lvert \Fix(\alpha) \rvert \leq
\frac{p^{\frac{n}{2}} \left(p^n-1\right)}{2}.
\]
\end{lemma}
\proof
Suppose $\alpha = \sigma^i D^j$.
We work in the cover $\SL_2(p^n)$.
Fixed points of $\alpha$ correspond to solutions of the system
\[
\begin{pmatrix}
a^{p^i} & b^{p^i} \\
c^{p^i} & d^{p^i}
\end{pmatrix}
=
\pm
\begin{pmatrix}
a & b \omega^j \\
c \omega^{-j} & d
\end{pmatrix}
\]
for $a,b,c,d \in \FF_{p^n}$.

If $i = 0$, then a fixed point can only be a diagonal matrix, and the number
of these in $\PSL_2(p^n)$ is at most $(p^n - 1) / 2$. This also covers the
case when $n = 1$.

Now assume that $i > 0$. Note that since $n$ is a power of $2$, $-1$ is a square
in $\FF_{p^n}$.
The diagonal fixed points of $\alpha$, i.e.,  the cases when $b = c = 0$,
correspond to fixed points of $\sigma^i$.
On the other hand, as long as the matrix is not diagonal, i.e.,  $bc \neq 0$, a
solution is possible if and only if $j$ is divisible by $p^i - 1$. Write $j = (p^i -1)k$ for some $k$, so that $\alpha = \sigma^i D^{(p^i -1)k}
= D^k \sigma^i D^{-k}$.
In this situation, an element $x \in \PSL_2(p^n)$ is fixed under $\alpha$
if and only if $x^{D^k}$ is fixed under $\sigma^i$. All in all, we therefore
have that
\[
\lvert \Fix(\alpha) \rvert =
\lvert \Fix(\sigma^i) \rvert \leq
\lvert \Fix(\sigma^{\frac{n}{2}}) \rvert =
\lvert \PSL_2(p^{\frac{n}{2}}) \rvert. \qedhere
\]
\endproof

\begin{lemma} \label{l:size_of_ad_map}
Let $\alpha \in \Aut(S)$. Set
\[
\ad_{\alpha} \colon S \to S, \qquad a \mapsto [a, \alpha].
\]
Then
\[
\lvert
\im \ad_{\alpha}
\rvert =
\frac{|S|}{\lvert \Fix(\alpha) \rvert}.
\]
\end{lemma}
\proof
For elements $a, b \in S$, we have
$\ad_{\alpha}(a) = \ad_{\alpha}(b)$ if and only if
$a b^{-1} = (a b^{-1})^\alpha$, which is the same
as saying that $a$ and $b$ belong to the same coset
of $\Fix(\alpha)$. The claim follows immediately.
\endproof

\subsection{$2$-Engel word: Variations}

Let $w = [X_1,X_2,X_2]$ be the $2$-Engel word. In expanded form, this is
\[
w = X_2^{-1} X_1^{-1} X_2 X_1 X_2^{-1} X_1^{-1} X_2^{-1} X_1 X_2 X_2.
\]

If a group satisfies the word $w$, it must be nilpotent. Therefore the word map of
$w$ on $\PSL_2(2) \cong S_3$ and $\Sz(2) \cong C_5 \rtimes C_4$ is not constant.
It follows from Lemma \ref{l:bor_reduction_simple_groups} that in order to
verify that $w$ is VSMB, we only need to consider the varied coset word maps
of $w$ on $\PSL_2(p^n)$ with $p > 2$ and $n$ a power of $2$.

It follows from Lemma \ref{l:bors_reductions_for_VSMB} that every variation of
$w$ in which appearances of $X_1^{\pm 1}$ are replaced by using more than one
second index are VSMB. Therefore we only need to vary occurrences of the variable
$X_2$. A general variation of $w$ can therefore be assumed to be of the form
\[
\tilde w = Y_1^{-1} X^{-1} Y_2 X Y_3^{-1} X^{-1} Y_4^{-1} X Y_5 Y_6
\]
with some of the $Y_i$ being potentially equal. Using Lemma \ref{l:bors_reductions_for_VSMB},
we can further reduce the words that need to be checked.
As long as
$|\{ Y_1, \dots, Y_6 \}| \geq 3$,
the variation itself is VSMB.
Therefore we can assume that either we are dealing with the original word
or with a variation in which each $Y_i$ is equal to a variable
$Z_1$ or $Z_2$.

For any such word $\tilde w$, let $S = \PSL_2(p^n)$ with $p > 2$ and $n$ a power
of $2$; we only need to inspect these simple groups by Lemma \ref{l:bor_reduction_simple_groups}.
Consider the coset word map induced by
elements $x, y_1, \dots, y_6$ in $\Aut(S)$.
Here we can assume that $x, y_i$ belong to the subgroup of $\Aut(S)$
generated by the field and diagonal automorphisms.
We want to show that this word map
is not constant. To this end, assume the contrary. For any $a, b_1, \dots, b_6
\in S$, we therefore have
\begin{equation} \label{e:coset_identity}
\tilde w(ax, b_1 y_1, \dots, b_6 y_6) = \tilde w(x, y_1, \dots, y_k).
\end{equation}

\subsection{$2$-Engel word: Inspecting the original}

Let us first deal with the original $2$-Engel word.
Insert $b_i = 1$ into \eqref{e:coset_identity}
and collect the left hand side to get that
\begin{equation} \label{e:2engel_collection_a}
a^{- x y}
a^{y^{-1} x y}
a^{- x y x^{-1} y^{-1} x y}
a^{y x y x^{-1} y^{-1} x y}
\equiv 1
\end{equation}
for all $a \in S$. Cancelling $x y$ and collecting, we obtain
\[
[a, y^{-1}] [a, y]^{x y x^{-1} y^{-1}} \equiv 1,
\]
and by additionally cancelling $y^{-1}$, it follows that
\[
[a, y]^{-1} [a, y]^{x y x^{-1}} \equiv 1.
\]
This can be rewritten as
\begin{equation} \label{e:2engel_collection_a_simplified}
[a, y, x y x^{-1}] \equiv 1,
\end{equation}
which is the same as saying that
\[
\im \ad_{y} \subseteq \Fix(x y x^{-1}).
\]
Comparing the sizes, it now follows from Lemma \ref{l:size_of_ad_map}
and Lemma \ref{l:out_fixed_points} that, as long as
$y$ is nontrivial,
\[
\frac{p^n \left( p^{2n} - 1 \right)}{2} =
\lvert \PSL_2(p^n) \rvert \leq
\lvert \Fix(y) \rvert \cdot \lvert \Fix(x y x^{-1}) \rvert \leq
\frac{p^n \left( p^n - 1 \right)^2}{4}
\]
This is impossible. Hence $y = 1$ and in this case
it is clear that \eqref{e:2engel_collection_a} can not hold.

\subsection{$2$-Engel word: Inspecting the variations}

In each of the proper variations $\tilde w$, we have a variable $X$
and two other variables $Z_1, Z_2$, each one with multiplicity $3$.
Consider \eqref{e:coset_identity} with the variables $ax, by, cz$.
Insert $a = c = 1$. After collecting, we obtain
\begin{equation} \label{e:collecting_the_engel_variation}
b^{\pm \alpha} b^{\pm \beta} b^{\pm \gamma} \equiv 1
\end{equation}
for some fixed $\alpha, \beta, \gamma \in \langle D, \sigma \rangle  \leq \Aut(S)$
depending on $x,y,z$. Insert the diagonal matrix $D^2 \in \PSL_2(p^n)$ into the last equality.
The automorphisms $\alpha, \beta, \gamma$ act on it as powers of the
Frobenius automorphism $\sigma$. We obtain
\[
(D^2)^{\pm p^i \pm p^j \pm p^k} = 1
\]
for some fixed $0 \leq i,j,k < n$. Note, however, that $D^2$ is of order $(p^n - 1)/2$.
This number is even as long as $n$ is a proper power of $2$. On the other hand, the sum
$\pm p^i \pm p^j \pm p^k$ is always odd. Thus we are forced into the conclusion $n = 1$.
By Lemma \ref{l:bor_reduction_simple_groups}, it suffices
to consider primes $p \leq 3$, as the maximum multiplicity
of a variable in $\tilde w$ is $4$.
Therefore it suffices to verify that $\tilde w$ can not satisfy
an identity of the form \eqref{e:coset_identity} in the group
$\PSL_2(3) \cong \Alt(4)$.
In this group, insert the element $(1 \ 2)(3 \ 4)$ into
\eqref{e:collecting_the_engel_variation}. This element is fixed by
$D$, which is represented as conjugation by $(1 \ 2)$. The automorphisms
$\alpha, \beta, \gamma$ act as powers of $D$, and hence we obtain
\[
\left( (1 \ 2)(3 \ 4) \right)^{\pm 1 \pm 1 \pm 1} = 1.
\]
a contradiction which completes our analysis of the variations.

\subsection{Metabelian word}

In this section, we deal with the metabelian word $[[X_1,X_2], [X_3, X_4]]$.
In expanded form, this is
\[
X_2^{-1} X_1^{-1} X_2 X_1
X_4^{-1} X_3^{-1} X_4 X_3
X_1^{-1} X_2^{-1} X_1 X_2
X_3^{-1} X_4^{-1} X_3 X_4.
\]
By Lemma \ref{l:bors_reductions_for_VSMB}, any proper variation of the metabelian
word is VSMB.
Therefore it suffices to consider only the original word.

Let $S$ be a nonabelian finite simple group.
Consider the coset word map induced by elements $x,y,z,t \in \Aut(S)$.
Assume that this map is constant on $S$. For any $a,b,c,d \in S$, we therefore
have
\[
[[ax,by], [cz,dt]] = [[x,y], [z,t]] = {\rm const}.
\]
Expand the first commutator to get
\[
[[a,by]^x [x,by], [cz, dt]] = {\rm const}.
\]
and once again
\[
[[a,by]^x, [cz, dt]]^{[x, by]} [[x,by], [cz, dt]] = {\rm const}.
\]
The second factor is constant, since it equal to the original
word map with $a = 1$. We conclude that
\[
[[a,by]^x, [cz, dt]]^{[x, by]} = {\rm const}.
\]
Inserting $a = 1$, we see that the value of the last word must in fact be equal
to $1$. Therefore
\begin{equation} \label{e:temp_derived_word}
[[a,by]^x, [cz, dt]] \equiv 1.
\end{equation}
Expand the first commutator, now in the second variable, to obtain
\[
[[a,y]^x[a,b]^{yx}, [cz, dt]] \equiv 1,
\]
and once again
\[
[[a,y]^x, [cz, dt]]^{[a,b]^{yx}} [[a,b]^{yx}, [cz, dt]] \equiv 1.
\]
We see that the first factor is trivial by inserting $b = 1$ into
\eqref{e:temp_derived_word}. It follows that
\begin{equation} \label{e:metabelian_expansion}
[[a,b]^{yx}, [cz, dt]] \equiv 1.
\end{equation}
Since $S$ is a perfect group, we now conclude that
the element $[cz, dt]$ must fix the whole of $S$, and so
$[cz, dt] = 1$. Therefore $[z,t] = [cz,t] = 1$,
forcing $[c,t] = 1$ and similarly $[z,d] = 1$.
This gives $z = t = 1$. A symmetric argument shows that
$x = y = 1$. But now $S$ should satisfy the metabelian identity,
a contradiction.

\subsection{Coset probabilities}
\label{subs:coset_probabilities}

For a group $G$, elements $g_1, g_2, \dots, g_d$ and a normal subgroup
$N$ of $G$, denote
\[
\Prob_{w = 1}^{(g_1, g_2, \dots, g_d)}(N)
=
\frac
{|\{ (n_1, n_2, \dots, n_d) \in N^d \mid
w(n_1 g_1, n_2 g_2, \dots, n_d g_d) = 1 \}|}
{|N|^d}.
\]
This is the coset probability of $N$ in $G$ of satisfying the word $w$.
Taking $N = G$, we recover the ordinary probability of satisfying $w$ in $G$.
As long as $N$ is not solvable, it is possible to
universally bound this coset probability.

\begin{proposition} \label{p:coset_gap}

Let $w$ be a VSMB word.
There exists a constant $\mu < 1$ such that whenever
$G$ is a finite group and $N$ its nonsolvable normal subgroup,
we have, for all $g_1, \dots, g_d \in G$,
\[
\Prob_{w = 1}^{(g_1, g_2, \dots, g_d)}(N)
\leq \mu \cdot \mathbbm{1}_{w(g_1, g_2, \dots, g_d) \in N}.
\]
\end{proposition}
\proof

As long as $\Prob_{w=1}^{(g_1, g_2, \dots, g_d)}(N) > 0$,
we must have $w(g_1, g_2, \dots, g_d) \in N$. This is the reason
why we include the indicator function in the statement.

Let $T^k$ be a chief factor of $N$ with $T$ a nonabelian simple group.
After replacing $G$ by a suitable quotient, we can assume that $C_G(T^k) = 1$,
and so $T^k \leq G \leq \Aut(T^k)$.

There exists $C > 0$ depending only on $w$ such that whenever $|T| \geq C$, we
have (see \cite[Theorem 4.5]{LarSha17} together with \cite[Lemma 2.7]{Bor17})
\[
\Prob_{w = 1}^{(g_1, g_2, \dots, g_d)}(N)
\leq \frac12.
\]
Therefore it suffices to consider only finitely many options for the simple
group $T$.

As $w$ is assumed to be a VSMB word, no varied coset word map on $T$ is
constant (see \cite[Definition 2.8]{Bor17}). There now exists an
$0 < \epsilon < 1$ depending only on $w$ and $T$ such that
(see \cite[Lemma 2.12]{Bor17})
\[
|\{ (n_1, n_2, \dots, n_d) \in N^d \mid
w(n_1 g_1, n_2 g_2, \dots, n_d g_d) = 1 \}|
\leq
\epsilon^{\lceil k / \ell(w)^2 \rceil} |T|^{kd}.
\]
Therefore there exists a $\bar C > 0$ such that whenever $k \geq \bar C$, we have
\[
\Prob_{w = 1}^{(g_1, g_2, \dots, g_d)}(N)
\leq
\frac12.
\]
Therefore it suffices to consider only finitely many options for the
multiplicity $k$ of $T$.

Now, as $G \leq \Aut(T^k)$, there are only finitely many options left for
the group $G$. None of these groups
satisfy a coset identity since $w$ is VSMB, and so each value
$
\Prob_{w = 1}^{(g_1, g_2, \dots, g_d)}(N)
$
is smaller than $1$. Thus we can take $\delta$ to
the maximum of all these values and $\frac12$.
\endproof

A consequence of the existence of this bound is the following
bound for the probability of satisfying $w$ when extending groups.

\begin{corollary} \label{c:probability_of_nonsolvable_extensions}

Let $w$ be a VSMB word.
There exists a constant $\mu < 1$ such that whenever
$G$ is a finite group and $N$ its nonsolvable normal subgroup,
we have
\[
\Prob_{w=1}(G) \leq \mu \cdot \Prob_{w=1}(G/N).
\]
\end{corollary}
\proof

Let $\mathcal{R}$ be a set of coset representatives of $N$ in $G$.
We have
\[
\Prob_{w=1}(G) =
\frac{1}{|G|^d}
\sum_{(r_1, r_2, \dots, r_d) \in \mathcal{R}^d}
|N|^d \cdot \Prob_{w=1}^{(r_1, r_2, \dots, r_d)}(N).
\]
Bounding the latter probability using Proposition \ref{p:coset_gap},
we obtain
\[
\Prob_{w=1}(G)
\leq
\frac{1}{|G/N|^d}
\sum_{(r_1, r_2, \dots, r_d) \in \mathcal{R}^d}
\mu \cdot \mathbbm{1}_{w(r_1, r_2, \dots, r_d) \in N}
=
\mu \cdot \Prob_{w=1}(G/N),
\]
as claimed.
\endproof

\begin{corollary} \label{c:chain}

Let $w$ be a VSMB word.
Let $G$ be a group with a chain
$N_1 > N_2 > \dots$
of normal subgroups of finite index in $G$
such that the consecutive factors $N_i / N_{i+1}$
are not solvable. Then
\[
\lim_{i \to \infty} \Prob_{w=1}(G/N_i) = 0.
\]
\proof
Immediate by Corollary \ref{c:probability_of_nonsolvable_extensions}.
\endproof

\end{corollary}

\section{Solvable groups with non-trivial action of the group (resp. derived subgroup)}
\label{s:solvable_groups}

In this section, we deal with bounding the probability $\Prob_{w=1}(G)$ for finite
solvable groups $G$. As explained above, this is reduced to bounding
$\Prob_{w=1}(G)$ in the case when the verbal subgroup of $w$ in $G$,
denoted throughout by $V$, is a minimal normal subgroup that
is a vector space over a finite field $\FF_p$,
say of dimension $n$. We will exploit this action, so we assume throughout
this section that $n > 1$. Our assumption in this section will be
that in the $2$-Engel case, the action of $G$ on $V$ is non-trivial,
and that in the metabelian case, the action of $G'$ on $V$ is non-trivial.

\subsection{General principle for bounding the probability}
\label{subs:general_principle_for_bounding_probability_higher_dimensional_verbal}

Let $\mathcal{R}$ be a set of coset representatives for $V$ in $G$.
The probability of satisfying $w$ in $G$ can be expressed as
\[
\Prob_{w=1}(G) =
\frac{1}{|G|^d}
\sum _{a_i \in V, r_i \in \mathcal{R}} \mathbbm{1}_{w(a_1 r_1, \dots, a_d r_d) = 1}.
\]
Each summand can be expanded as
\[
w(a_1 r_1, \dots, a_d r_d) =
\prod_{i = 1}^d a_i^{w_i(r_1, \dots, r_d)} \cdot
w(r_1, \dots, r_d)
\]
for some endomorphisms $w_i(r_1, \dots, r_d) \in \End(V)$.
Set
\[
\texttt{BAD} =
\{ (r_1, \dots, r_d) \in \mathcal{R}
\mid
\forall i. \; w_i(r_1, \dots, r_d) = 0_{\End(V)} \}.
\]
This set consists of those tuples of elements of $\mathcal{R}$ for which a summand above
is independent of the values $a_i \in V$. Thus, these tuples are providing a coset
identity. Correspondingly, set $\texttt{GOOD} = \mathcal{R}^d - \texttt{BAD}$.
By first summing over the bad representatives, we have
\[
\frac{1}{|G|^d}
\sum_{(r_1, \dots, r_d) \in \texttt{BAD}} |V|^d \cdot \mathbbm{1}_{w(r_1, \dots, r_d) = 1}
\leq \frac{|\texttt{BAD}|}{|\mathcal{R}|^d}.
\]
On the other hand, for a good tuple of representatives,
at least one exponential endomorphism,
say $w_j(r_1, \dots, r_d)$, acts nontrivially on $V$.
Its kernel in $V$ is therefore of codimension at least $1$.
In this case, we have
\[
| \{ a_j \mid
w(a_1 r_1, \dots, a_d r_d) = 1
\} |
\leq
|C_V(w_j(r_1, \dots, r_d))|
\leq \frac{|V|}{p},
\]
and it follows form this that by summing over the good representatives, we have
\[
\frac{1}{|G|^d}
\sum_{(r_1, \dots, r_d) \in \texttt{GOOD}}
\sum_{a_i \in V}
\mathbbm{1}_{w(a_1 r_1, \dots, a_d r_d) = 1}
\leq
\frac{1}{|G|^d}
\sum_{(r_1, \dots, r_d) \in \texttt{GOOD}}
|V|^{d-1} \frac{|V|}{p}
=
\frac{|\texttt{GOOD}|}{p|\mathcal{R}|^d}
\]
We can collect the two upper bounds to finally obtain
\[
\Prob_{w=1}(G) \leq
\frac{|\texttt{BAD}|}{|\mathcal{R}|^d}
+
\frac{|\texttt{GOOD}|}{p |\mathcal{R}|^d}.
\]
Taking $|\texttt{BAD}| + |\texttt{GOOD}| = |\mathcal{R}|^d$ into account,
we can take the latter one step further and write
\[
\Prob_{w=1}(G) \leq
\frac{1}{p} +
\left( 1 - \frac{1}{p} \right) \frac{|\texttt{BAD}|}{|\mathcal{R}|^d}
\leq
\frac12 \left( 1 + \frac{|\texttt{BAD}|}{|\mathcal{R}|^d} \right).
\]
We will use this general principle for bounding the word probability.
In order for it to give us a gap on word probability, we will need to show
that for a given word $w$, there is a gap on the relative size of the set
$\texttt{BAD}$ inside $\mathcal{R}^d$.

\subsection{2-Engel word: Inspecting badness}

In this section, we focus on the case of the $2$-Engel word $[ax, by, by]$
with $a,b \in V$ and $x,y \in \mathcal{R}$. In order to obtain the equations
for defining the $\texttt{BAD}$ representatives, 
substitute $b := 1$ (resp. $a := 1$) and collect the resulting expression.
We first get, as in \eqref{e:2engel_collection_a} and simplified to
\eqref{e:2engel_collection_a_simplified},
\[
[a, y, x y x^{-1}] \equiv 1,
\]
so that the operator $(1 - y)(1 - x y x^{-1})$ acts trivially on $V$.
On the other hand, the condition that $[x,by,by]$ be constant
can be translated by expanding commutators into
\[
[x,y,b]^y [[x,b]^y, y] \equiv 1.
\]
Collecting the exponents at $b$, we see that the operator
$(1 - [x,y]) - (1 - x ) (1 - y)$ must also annihilate
everything on $V$. Now, since $[x,y]$ commutes with $y$
in its action on $V$, it follows that
$(1-x)(1-y)$ also commutes with $y$,
and hence $x(1-y)$ also commutes with $y$.
Thus we obtain
\[
0 \equiv (1-y)(1 - x y x^{-1}) =
(1-y)x (1-y) x^{-1} = (1 - y)^2.
\]
The latter means that $y$ acts quadratically on $V$, i.e.,
for all $a \in V$ we have $[a,y,y] = 1$.
Therefore
\[
[a, y^p] = [a, y]^p = 1,
\]
and so $y^p \in C_G(V)$. We have thus derived the inclusion
\[
\texttt{BAD} \subseteq
\mathcal{R} \times \{ y \in \mathcal{R} \mid y^p \in C_G(V) \}.
\]

\subsection{$2$-Engel word: Bounding badness in a nontrivial action}
\label{subs:2engel_bounding_badness}

Suppose that the action of $G$ on $V$ is nontrivial, that is $C_G(V) \neq G$.
We can view $V$ as a modular irreducible representation of $G / C_G(V)$. Since
$G / C_G(V)$ satisfies the $2$-Engel word, it is nilpotent. Now, $G / C_G(V)$
can not be a $p$-group, since the only irreducible representation of a
$p$-group in characteristic $p$ is the trivial one. Hence
the Sylow $p$-subgroup $P/C_G(V) \leq G/C_G(V)$ is proper.
We can identify $\mathcal{R}$ with cosets of $V$ in $G$,
and in this sense
\[
\texttt{BAD} \subseteq
\mathcal{R} \times (P/V).
\]
This gives the desired bound for the relative size of bad representatives,
\[
\frac{|\texttt{BAD}|}{|\mathcal{R}|^2}
\leq
\frac{|P|}{|V| |\mathcal{R}|} =
\frac{1}{|G : P|} \leq \frac12.
\]

\subsection{Metabelian word: Inspecting badness}

In this section, we focus on the case of the metabelian word
$[[ax, by], [cz, dt]]$ with $a,b,c,d \in V$ and $x,y,z,t \in \mathcal{R}$.
Collecting each word value separately, we obtain
\[
(1-x) y (1 - [z,t]) \equiv (1-y) x (1 - [z,t]) \equiv 0
\]
and similarly for the symmetric situation. Set $A = 1 - [z,t]$
and $B = 1 - [x,y]$. Thus $xyA = yA$ and $yxA = xA$.
Note that $yx - xy = yx B$. Now, since
$[x,y]$ and $[z,t]$ induce commuting operators on $V$, it follows
that
\[
(x-y)A = (yx - xy)A = yx B A = yx A B = x A B.
\]
This gives
\[
yA = xA - xAB = x A [x,y] = x[x,y] A = y^{-1} x y A = y^{-1} y A = A.
\]
(Similarly we can derive other equalities.) This means that
we have, for all $a \in V$,
\[
[[a, y], [z,t]] = 1.
\]
Thus we have the inclusion
\[
\texttt{BAD} \subseteq
\mathcal{R} \times
\{ (y,z,t) \in \mathcal{R}^3 \mid [V, y, [z,t]] = 1 \}.
\]

\subsection{Metabelian word: Bounding badness in a non-trivial action
of derived subgroup}

Suppose that the action of $G'$ on $V$ is nontrivial, that is $G' \not\subseteq C_G(V)$.
We can identify $\mathcal{R}$ with cosets of $V$ in $G$.
Let
\[
\texttt{Y}_{[z,t]} = \{ y \in G/V \mid [V, y, [z,t]] = 1 \}.
\]
Note that for $y_1, y_2 \in \texttt{Y}_{[z,t]}$, we have
\[
[a, y_1 y_2, [z,t]] =
[a, y_2, [z,t]]^{[a, y_1]^{y_2}} \cdot
[a, y_1, [z,t]]^{[a, y_1, y_2]} \cdot
[[a, y_1], y_2, [z,t]] =
1,
\]
and so $\texttt{Y}_{[z,t]}$ is a subgroup of $G/V$.
As long as $[z,t]$ is not trivial in $G / C_G(V)$,
this is a proper subgroup, since $[V,G] = V$.
Set
\[
\texttt{UGLY} = \{ (z,t) \in (G/V)^2 \mid [z,t] \in C_G(V) \}.
\]
Thus we can express
\[
\texttt{BAD} =
\left( \mathcal{R}^2 \times \texttt{UGLY} \right)
\cup
\left( \mathcal{R} \times
\bigcup_{(z,t) \in \mathcal{R}^2 - \texttt{UGLY}}
\texttt{Y}_{[z,t]}
\right)
\]
and compute, taking into account that $\texttt{Y}_{[z,t]}$ is of index
at least $2$ in $G/V$,
\[
|\texttt{BAD}|
\leq
|\mathcal{R}|^2 |\texttt{UGLY}|
+
|\mathcal{R}|
\left(
|\mathcal{R}|^2 - |\texttt{UGLY}|
\right)
\frac{|G / V|}{2}.
\]
Thus we obtain a bound for the badness ratio,
\[
\frac{|\texttt{BAD}|}{|\mathcal{R}|^4}
\leq
\frac12 + \frac12 \frac{|\texttt{UGLY}|}{|\mathcal{R}|^2}.
\]
Since $G / C_G(V)$ is assumed not to be abelian, we obtain the last bound
\[
\frac{|\texttt{UGLY}|}{|\mathcal{R}|^2} =
\frac{|\{ (z,t) \in (G/C_G(V))^2 \mid [z,t] = 1 \}|}{|G / C_G(V)|^2}
\leq \frac58.
\]

\section{Solvable groups with trivial action of the group (resp. derived subgroup)}
\label{s:1dim_verbal}

In this case, every tuple in $\mathcal{R}^4$ is bad, and so our general
procedure for bounding the word probability in terms of counting bad tuples
does not work. We will therefore use the following principle.

\subsection{General principle for bounding the probability}
\label{subs:general_principle_for_bound_probability_1_dimension}

The probability of satisfying $w$ in $G$ can be expressed as
\[
\Prob_{w=1}(G) =
\frac{1}{|G|^d}
\sum _{g_2, \dots, g_d \in G} |\{ g_1 \in G \mid w(g_1, \dots, g_d) = 1 \}|.
\]
Let $\texttt{BAD} \subseteq G^{d-1}$ be a certain subset of tuples,
and set correspondingly $\texttt{GOOD} = G^{d-1} - \texttt{BAD}$.
Denote
\[
C_w(g_2, \dots, g_d) = \{ g_1 \in G \mid w(g_1, \dots, g_d) = 1 \}.
\]

\noindent {\bf Assumption:}
There exist absolute constants
$0 < \delta_{\texttt{GOOD}}, \delta_{\texttt{BAD}} < 1$,
depending only on $w$ and not on $G$,
such that:
\begin{enumerate}
\item
$\forall (g_2, \dots, g_d) \in \texttt{GOOD}. \;
|C_w(g_2, \dots, g_d)| \leq \delta_{\texttt{GOOD}} \cdot |G|$

\item
$|\texttt{BAD}| \leq \delta_{\texttt{BAD}} \cdot |G|^{d-1}$
\end{enumerate}

Under the above assumption, we can bound the word
probability as follows. First of all, we let the sum expressing
the word probability go over the good and the bad tuples
separately,
\begin{align*}
\Prob_{w = 1}(G)
& \leq
\frac{1}{|G|^d}
\sum _{g_2, \dots, g_d \in \texttt{GOOD}} |C_w(g_2, \dots, g_d)|
+
\frac{1}{|G|^d}
\sum _{g_2, \dots, g_d \in \texttt{BAD}} |G| \\
& \leq
\frac{1}{|G|^d} |\texttt{GOOD}| \delta_{\texttt{GOOD}} |G|
+
\frac{1}{|G|^d} |\texttt{BAD}| |G|.
\end{align*}
Taking $|\texttt{BAD}| + |\texttt{GOOD}| = |G|^{d-1}$ into account,
we can therefore write
\[
\Prob_{w=1}(G) \leq
\delta_{\texttt{GOOD}} +
( 1 - \delta_{\texttt{GOOD}}) \frac{|\texttt{BAD}|}{|G|^{d-1}}.
\]
we can bound the relative badness by assumption, and hence
\[
\Prob_{w=1}(G) \leq
\delta_{\texttt{GOOD}} +
( 1 - \delta_{\texttt{GOOD}}) \delta_{\texttt{BAD}},
\]
which gives an absolute upper bound on the word probability in $G$.

Our method of satisfying the above assumption will be the following. In order
to obtain the constant $\delta_{\texttt{BAD}}$,  we will repeatedly use the
fact that a proper subgroup of a group is of index at least $2$. This is the
explanation for the gap in this case. As for the constant 
$\delta_{\texttt{GOOD}}$, 
we will consider the word map $w(\cdot, g_2, \dots, g_d) \colon G \to G$
and bound its fiber over $1$ using the following.

\begin{lemma}[Fiber of restricted homomorphism] \label{l:fiber_of_restricted_homomorphism}
Let $\phi \colon G \to G$ be a map. Suppose there exists a subgroup
$H \leq G$ for which we have
\[
\forall g \in G \ \forall h \in H. \; \phi(gh) = \phi(g) \phi(h)
\quad
\text{and}
\quad
\phi(H) \neq \{ 1 \}.
\]
Then
\[
|\phi^{-1}(1)| \leq \frac{1}{2} |G|.
\]
\end{lemma}
\proof
Let $R$ be a set of coset representatives of $H$ in $G$. We can express
\[
|\phi^{-1}(1)| = \sum_{r \in R} |\{ h \in H \mid \phi(rh) = 1 \}|.
\]
By our assumption on $\phi$, the condition $\phi(rh) = 1$ is 
equivalent to $\phi(h) = \phi(r)^{-1}$. Since the restriction
$\phi \lvert_H \colon H \to G$ is a homomorphism,
the size of any fiber is at most the size of the kernel. Hence
we have
\[
|\phi^{-1}(1)| \leq \sum_{r \in R} |\ker \phi \lvert_H|.
\]
Now, as $\phi \lvert_H$ is assumed to be non-trivial, 
its kernel is of index at least $2$ in $H$. It now follows that
\[
|\phi^{-1}(1)| \leq |R| \cdot \frac12 |H| = \frac12 |G|. \qedhere
\]
\endproof

\subsection{Long commutator}

Consider the long commutator word,
\[
\gamma_d(X_1, \dots, X_d) = [X_1, \gamma_{d-1}(X_2, \dots, X_d)]
= [X_1, [X_2, [\dots, X_d]]].
\]
We know that in a nonabelian group, the fiber over $1$
of $\gamma_2(X_1, X_2)$ is of relative size at most $\frac58$.
Similarly, it is known (see \cite{Erf07}) that
there exists a probability gap for the long commutator word.
We give a sample application of our general principle for bounding the probability by providing an alternative proof of this fact.

Let us show by induction that a bound exists for all the long
commutator words.
Let $G$ be a group that does not satisfy the word $\gamma_d$.
Set
\begin{align*}
\texttt{BAD}
&= \{ (g_2, \dots, g_d) \in G^d \mid
C_w(g_2, \dots, g_d) = G
\} \\
&= \{ (g_2, \dots, g_d) \in G^d \mid
\gamma_{d-1}(g_2, \dots, g_d) \in Z(G)
\}.
\end{align*}
The size of the latter set is
\[
|\texttt{BAD}| =
|\{ (g_2, \dots, g_d) \in (G/Z(G))^d \mid
\gamma_{d-1}(g_2, \dots, g_d) = 1 \}|
\cdot |Z(G)|.
\]
The group $G / Z(G)$ does not satisfy the word $\gamma_{d-1}$.
Therefore we can argue by induction that there is a constant
$\delta_{d-1}$ with
\[
|\{ (g_2, \dots, g_d) \in (G/Z(G))^{d-1} \mid
\gamma_{d-1}(g_2, \dots, g_d) = 1 \}| \leq
\delta_{d-1} |G/Z(G)|^{d-1}.
\]
Thus
\[
|\texttt{BAD}| \leq \delta_{d-1} |G|^{d-1},
\]
so we can take $\delta_{\texttt{BAD}} = \delta_{d-1}$.
On the other hand, for a tuple $(g_2, \dots, g_d) \notin \texttt{BAD}$,
we have
\[
C_w(g_2, \dots, g_d) = C_G(\gamma_{d-1}(g_2, \dots, g_d)),
\]
which is a proper subgroup of $G$, and therefore
\[
|C_w(g_2, \dots, g_d)| \leq \frac12 |G|.
\]
Therefore we can take $\delta_{\texttt{GOOD}} = \frac12$.
This gives a bound for the probability,
\[
\Prob_{\gamma_d = 1}(G) \leq \frac12 + \frac12 \delta_{d-1} =: \delta_d.
\]
Since $\delta_2 = \frac58$, we can derive inductively that
we can take $\delta_d = 1 - \frac{3}{2^{d+1}}$.
In this case, the obtained bound is sharp, as can be seen for example by looking
at dihedral groups.

\subsection{$2$-Engel word: Bounding two types of badness}
\label{subs:engel_bounding_two_types}

We can assume that $V = \FF_p$ and that $V$ is central in $G$.  The case when
$G$ acts nontrivially on $V$ has been dealt with in Section
\ref{subs:2engel_bounding_badness}. Since $G/V$ is assumed to be $2$-Engel, it is
nilpotent, and so $G$ must also be nilpotent. As $V$ is the smallest normal
subgroup of $G$, this implies that $G$ must in fact be a $p$-group. Moreover,
$G / V$ is of nilpotency class at most $3$, so $G$ is of nilpotency class at most
$4$. We will use this fact freely in what follows. 
Moreover, it is known that $2$-Engel $p$-groups with $p \neq 3$ are actually
nilpotent of class at most $2$ (not just $3$) since they always satisfy the identity $[X,Y,Z]^3 = 1$. For the particular case when $p = 3$, we will
make use of the following.

\begin{claim}
If $p = 3$, then
\[
\forall y \in G. \; [G, y, G, y] = 1
\quad
\Longleftrightarrow
\quad
\forall y \in G. \; [G', y, y] = 1,
\]
and these two equivalent conditions imply that
$[G,G,G] \leq V$.
\end{claim}
\proof
We know that $[G,G,G,G]^3 = 1$,
and that $[G,G,G]^3 \leq V$.
So we have (see \cite[Lemma 2.2(v)]{GupNew89})
\[
[a,y,b,y] = [a,b,y,y]^{-1}.
\]
for all $a,b \in G$. This proves the equivalence in the claim.
As for the second part of the claim, suppose that
$\gamma_4(G) \neq 1$. Thus $\gamma_4(G) = V$.
For any $c \in G$, we have
\[
1 \equiv [c, ab, ab] \equiv [c, a, b] [c, b, a] \pmod{\gamma_4(G)},
\]
and using the Jacobi identity modulo $\gamma_4(G)$, we obtain
\[
[c, a, b] \equiv [a, b, c] \pmod{\gamma_4(G)}.
\]
Thus we have cyclic invariance of commutators of length $3$.
Now we have, by \cite[Vol. 2, p. 43]{Rob72},
that for any $d \in G$,
\[
[[a,b],c,d]=[[a,b],d,c]^{-1}
\quad
\text{and}
\quad
[[a,b],[c,d]]=[[a,b],c,d]^2.
\]
Using all the above, we can execute the computation
\begin{align*}
[a,b,c,d]^2
&= [a,b, [c,d]] \\
&= [[c,d], [a,b]]^{-1} \\
&= [c,d,a,b]^{-2} \\
&= [c,d,a,b] \\
&= [a,c,d,b] \\
&= [a,c,b,d]^{-1} \\
&= [c,a,b,d] \\
&= [a,b,c,d],
\end{align*}
giving $[a,b,c,d] = 1$. Hence indeed $\gamma_4(G) = 1$.
As $\gamma_3(G)$ is central in $G$, it must be cyclic,
since $V$ is the smallest normal subgroup of $G$.
But, since $\gamma_3(G)$ is of exponent $3$,
it follows that $\gamma_3(G) = V$, and the proof is complete.
\endproof

Let $y \in G$ and consider the map
\[
\phi_y \colon G \to V,
\qquad
a \mapsto [a,y,y].
\]
This may not be a homomorphism, but it does satisfy the following
expansion law which will be of use:
\begin{equation} \label{e:2engel_expansion_of_phi}
\phi_{y}(ab) =
\phi_y(a) \cdot
\phi_y(b) \cdot
[a,y,b,y].
\end{equation}
We will be interested in two possible situations, depending on
whether or not $\phi_y$ is a homomorphism.

\subsubsection{The nice situation: $\forall y \in G. \; [G, y, G, y] = 1$}

In this case, $\phi_y$ is a homomorphism for all choices of $y \in G$.
Note that this means that $\phi_y$ factors through the Frattini
quotient $G / \Phi(G)$, and so we can think of $\phi_y$ as a linear
functional over $\FF_p$ mapping into $V = \langle z \rangle$.
Let $\{ g_1, g_2, \ldots, g_d \}$ be a minimal generating set of $G$.
Every element $a \in G$ has a unique expansion
\[
a \equiv \prod_{i = 1}^d g_i^{\beta_i} \pmod{\Phi(G)}
\]
with $\beta_i \in \FF_p$.
We can express the map
$\phi_y$ in terms of the generating set of $G$.
Note that $\phi_y$ does not depend on the specific coset
representative of $y$ modulo $\Phi(G)$.
Set $[g_i, g_j, g_k] = z^{\gamma_{ijk}}$ with $\gamma_{ijk} \in \FF_p$.
Then we have
\[
\phi_y(g_i) = z^{\sum_{j,k} \gamma_{ijk} \beta_j \beta_k}.
\]
Set
\[
\texttt{BAD} = \{ y \in G \mid \phi_y \equiv 1 \}.
\]
Therefore $y \in \texttt{BAD}$ if and only if each of the $d$ quadratic forms
$\sum_{j,k} \gamma_{ijk} \beta_j \beta_k$ vanishes.
As $G$ is assumed not to be $2$-Engel, at least one of
these forms is not identically equal to $0$.
This form can be diagonalized (see \cite{Ser78}) to a form
\[
\sum_{l=1}^{d' - 1} \beta_l^2 + a_{d'} \beta_{d'}^2
\]
for some $1 \leq d' \leq d$ and $a_{d'} \in \FF_p$. The number of zeros of this
form can be bounded from above as follows.  For each of $\beta_2,\ldots,\beta_d$
there are at most $p$ choices and, after fixing these, there are at most two
possibilities for $\beta_1$ if $p\neq 2$, and at most one choice for $\beta_1$
if $p=2$. This implies that, if $p\neq 2$,
\[
|\texttt{BAD}| \leq 2 p^{d-1} \lvert \Phi(G) \rvert = \frac{2 |G|}{p} \leq \frac{2 |G|}{3},
\]
and similarly, if $p=2$, then $|\texttt{BAD}| \leq |G|/2$.

\subsubsection{The other situation: $\exists y \in G. \; [G, y, G, y] \neq 1$}

Note that this case is only possible when $p = 3$ (see \cite[Theorem
7.15]{Rob72}).  It follows from the claim stated at the begining of Subsection \ref{subs:engel_bounding_two_types}
that the restriction of the map $\phi_y$ to $G'$ is not always trivial.
Set
\[
\texttt{BAD'} = \{ y \in G \mid \phi_y(G') \equiv 1 \}.
\]
The same argument as above with $G$ replaced by $G'$ gives
the bound
\[
|\texttt{BAD'}| \leq \frac{2|G|}{3}.
\]

\subsection{$2$-Engel word: Bounding good fibers}

In order to provide a bound for the word probability,  we now need to ensure
that as long as $y$ does not belong to $\texttt{BAD}$ (or $\texttt{BAD’}$), we
can bound the number of solutions of the equation
\[
[a,y,y] = 1
\]
for $a \in G$. This is equivalent to saying that we want to provide a relative
upper bound for the fiber $\phi_y^{-1}(1)$. Note that the map $\phi_y$ is not
trivial in this situation.  We will need to analyse two cases.

\subsubsection{The nice situation: $\forall y \in G. \; [G, y, G, y] = 1$}

In this case, let $y$ be an element of $G$ that is not in $\texttt{BAD}$.
Thus $\phi_{y}$ is a non-trivial homomorphism, and so 
we are done by the Fiber of restricted homomorphism lemma (with $H = G$).

\subsubsection{The other situation: $\exists y \in G. \; [G, y, G, y] \neq 1$}

In this case, we will bound the relative size of $\phi_y^{-1}(1)$ for $y$ that
is not in $\texttt{BAD'}$.

\begin{claim}
Let $g \in G$ and $h \in G'$. Then
\[
\phi_y(gh) = \phi_y(g) \phi_y(h).
\]
\end{claim}
\proof
Expanding based on \eqref{e:2engel_expansion_of_phi},
we obtain
\[
\phi_y(g \cdot h) =
\phi_y(g) \cdot
\phi_y(h) \cdot
[g,y,h,y].
\]
As $h \in G'$, the last commutator is trivial.
\endproof

By our assumption and the claim at the beginning of Subsection
\ref{subs:engel_bounding_two_types}, the restriction $\phi_y \lvert_{G'}$ is
a non-trivial homomorphism. Therefore its fiber over $1$
can be bounded by the Fiber of restricted homomorphism lemma 
(with $H = G'$).

\subsection{Metabelian word: Bounding two types of badness}

We are in the situation when $G' \leq C_G(V)$. This means that $V = G''$
commutes with $G'$, and so $G'$ is nilpotent of class at most $2$. Since $V$
is the smallest normal subgroup of $G$, this implies that $G'$ must be a
$p$-group.

Let $y,z,t \in G$ and consider the map
\[
\phi_{y,z,t} \colon G \to V,
\qquad
a \mapsto [[a, y], [z,t]].
\]
This may not be a homomorphism, but it does satisfy the following
expansion law which will be of use:
\begin{equation} \label{e:metabelian_expansion_of_phi}
\phi_{y,z,t}(ab) = \phi_{y,z,t}(a) \cdot \phi_{y,z,t}(b) \cdot [[a,y,b],[z,t]].
\end{equation}
Set
\[
\texttt{BAD} = \{ (y,z,t) \in G \mid \phi_{y,z,t} \equiv 1 \}.
\]
With a fixed pair $(z,t)$, set
\[
S_{z,t} = \{ y \in G \mid \phi_{y,z,t} \equiv 1 \}.
\]
Note that for any $a \in G$, we have
\[
\phi_{y_1 y_2, z, t}(a) =
[[a, y_1 y_2], [z,t]] =
\phi_{y_2, z, t}(a) \cdot
\phi_{y_1, z, t}(a) \cdot
\phi_{y_2, z, t}([a, y_1]),
\]
so that $S_{z,t}$ is a subgroup of $G$. Thus we either have that $S_{z,t} = G$
or it is a proper subgroup of $G$.
The first case occurs if and only if $[z,t]
\in C_G(G')$.
Set
\[
\texttt{UGLY} =
\{ (z,t) \in G^2 \mid [z,t] \in C_G(G') \}.
\]
Thus we can express
\[
\texttt{BAD} =
\left( G \times \texttt{UGLY} \right)
\cup
\{
(y,z,t) \in G^3 \mid (z,t) \notin \texttt{UGLY}, \ y \in S_{z,t}
\},
\]
and so we have
\[
|\texttt{BAD}| =
|G| |\texttt{UGLY}|
+
\sum_{(z,t) \notin \texttt{UGLY}} |S_{z,t}|.
\]
When $(z,t)$ not an ugly pair, $S_{z,t}$
is of index at least $2$ in $G$. Thus we can bound
\[
|\texttt{BAD}| \leq
|G| |\texttt{UGLY}|
+
(|G|^2 -|\texttt{UGLY}|) \frac{|G|}{2}
=
\frac12 |G|^3 + \frac12 |G| |\texttt{UGLY}|.
\]
Note that as $G$ is assumed not to be metabelian,
$G / C_G(G')$ is not abelian. This means that we
can bound the relative number of ugly pairs,
finally giving
\[
|\texttt{BAD}| \leq \left( \frac12 + \frac12 \cdot \frac 58 \right) |G|^3
= \frac{13}{16} |G|^3.
\]

Later on, we will provide a bound for the good fibers. To this end, we will need to distinguish two cases, and one of these will require us to deal with the particular situation
when $[G', G, G']$ is a nontrivial subgroup of $G$.
This means that $G'$ is not contained in $C_G([G', G])$,
and so $G / C_G([G', G])$ is not abelian.
In this case, we will need to resort to the sets
\[
\texttt{BAD'} =
\{ (y,z,t) \in G^3 \mid \phi_{y,z,t}(G') = 1 \}
\]
and
\[
\texttt{UGLY'} =
\{ (z,t) \in G^3 \mid [z,t] \in C_G([G', G]) \}.
\]
The same argument as above gives the same bound
\[
|\texttt{BAD'}| \leq \frac{13}{16} |G|^3.
\]

\subsection{Metabelian word: Bounding good fibers}

In order to provide a bound for the word probability, we now need
to ensure that as long as a tuple $(y,z,t)$ does not belong to
$\texttt{BAD}$ (or $\texttt{BAD'}$),
we can bound the number of solutions of the
equation
\[
[[a,y], [z,t]] = 1
\]
for $a \in G$. This is equivalent to saying that we want
to provide a relative upper bound for the fiber
$\phi_{y,z,t}^{-1}(1)$. Note that the map $\phi_{y,z,t}$
is not trivial in this situation. We will need to
analyse two cases.

\subsubsection{The nice situation: $[G', G, G'] = 1$.}

In this case, we take a tuple that is not in $\texttt{BAD}$.
Thus $\phi_{y,z,t}$ is a nontrivial homomorphism, and so 
we are done by the Fiber of restricted homomorphism lemma (with $H = G$).

\subsubsection{The other situation: $[G', G, G'] \neq 1$.}

In this situation, we need to bound the relative
fiber size $\phi_{y,z,t}^{-1}(1)$ for a tuple
$(y,z,t) \notin \texttt{BAD'}$.

\begin{claim}
Let $g \in G$ and $h \in G'$. Then
\[
\phi_{y,z,t}(gh) = \phi_{y,z,t}(g) \phi_{y,z,t}(h).
\]
\end{claim}
\proof
Expanding based on \eqref{e:metabelian_expansion_of_phi},
we obtain
\[
\phi_{y,z,t}(g \cdot h) =
\phi_{y,z,t}(g) \cdot
\phi_{y,z,t}(h) \cdot
[[g,y,h], [z,t]].
\]
As $h \in G'$, the last commutator is trivial.
\endproof

By our assumption, the restriction $\phi_{y,z,t} \lvert_{G'}$ is
a non-trivial homomorphism. Therefore its fiber over $1$
can be bounded by the Fiber of restricted homomorphism lemma 
(with $H = G'$).

\newpage

\end{document}